\documentclass{amsart}
        \newtheorem{theorem}{{\bf Theorem}}
        
        \newtheorem{lemma}[theorem]{{\bf Lemma}}
        \newtheorem{proposition}[theorem]{{\bf Proposition}}
        \newtheorem{corollary}[theorem]{{\bf Corollary}}
        \newtheorem{definition}[theorem]{{\bf Definition}} %please number definitions
\theoremstyle{definition}
        \newtheorem{remark}[theorem]{{\bf Remark}}
        
\usepackage{amsxtra}
\usepackage[all]{xy}
\date{\today} 
\begin{document}
\title{Ramifications, Reduction of Singularities, Separatrices} 
\author{P. Fortuny Ayuso} 
\address{Department of Mathematical Sciences. School of Mathematics.
  Queen Mary College, University of London. Mile End Road. E1 4NS London, UK.}
\email{p.fortuny@maths.qmul.ac.uk}
%\ead{pf@maths.qmul.ac.uk} 
\thanks{Supported by the European Commission Marie
  Curie Fellowship Programme.}
\begin{abstract}
  We study the behaviour of sequences of blowing-ups under ramifications, and use
  the results to give a simple proof of Camacho-Sad's Theorem on the existence of
  Separatrices for singularities of plane holomorphic foliations. The main result
  we prove is that for any finite sequence $\pi$ of blowing-ups, there is a ramification
  morphism $\rho$ such that the elimination of indeterminations $\tilde{\pi}$ of $\pi^{-1}\circ \rho$
  is a sequence of blowing-ups \emph{with centers at regular points of the exceptional
  divisors}; moreover, we show that if $\pi$ is the reduction of singularities of a
  foliation $\mathcal{F}$, then $\rho$ can be such that $\tilde{\pi}^{\star}\rho^{\star}\mathcal{F}$
  has only simple singularities.
\end{abstract}
\maketitle
\section{Introduction}\label{sec:introduction}
In 1982, C. Camacho and P. Sad \cite{Camacho-Sad} proved that any holomorphic foliation
in $({\mathbb C}^{2},0)$ has at least one invariant anaylitc curve passing
through $(0,0)$. This result was searched for, at least, since XIX Century
(cf. \cite{Briot-Bouquet} in which it was proved for what are now called
\emph{simple} singularities). Afterwards, the same Camacho generalized it to
vector fields in singular surfaces whose dual resolution graph is a tree
\cite{Camacho}. In 1993, J. Cano \cite{CanoJ2} gave a
constructive proof, which he simplified in 1997 \cite{CanoJ4}. Camacho and Sad's
requires a detailed (and rather cumbersome) study of the combinatorial behaviour
of what has since been called the \emph{Camacho-Sad index} of an invariant
manifold. J. Cano's, on the other hand, although being much simpler, introduces
a somewhat artificial condition on that index.\par
Our work relies on the ---new, as far as the author knows--- study of the
behaviour of sequences of blowing-ups under ramification morphisms. These maps
were used successfully in \cite{Fortuny-Sanz} regarding Thom's Gradient
Conjecture (see \cite{Kurdyka-Mostowski-Parusinski}) and \cite{Corral} for
studying some properties of the polar curves of a foliation. In the present
paper we show ---roughly speaking--- how a sequence of blowing-ups can be modified
by a ramification to obtain a new sequence \emph{without combinatorial
  terms}: in the same way as one can modify a (singular) plane curve by a ramification in
order to get a union of regular curves, we ``modify'' a sequence of blowing-ups
in order to get another one which is ``regular''.
\section{Reduction of Singularities and Ramifications}\label{sec:ramifications}
Recall that a foliation $\mathcal{F}$ in a smooth analytic surface $\mathcal{X}$
can be
identified, \emph{locally}, with a germ of holomorphic $1-$form $\omega $, whose
local expression $a(x,y)dx + b(x,y)dy$ is such that $a$ and $b$ are relatively
prime in ${\mathbb C}\{x,y\}$. We will (abusing notation) speak indifferently of
$\mathcal{F}$ and $\omega$ when there is no possibility of confusion. A
\emph{separatrix} of $\mathcal{F}$ is a (germ of) holomorphic curve $\gamma:({\mathbb C},0) \rightarrow \mathcal{X}$
such that the pull-back $\gamma ^{\star}\omega $ is null. A point $P$ is a
\emph{simple singularity} of $\omega $ if there are coordinates at $P$ such that
the linear part of the vector field $-b(x,y)\partial /\partial x+a(x,y)\partial /\partial y$ has
two different eigenvalues $\mu \neq\lambda \neq 0$ with $\mu /\lambda \not\in
{\mathbb Q}_{>0}$.
\begin{definition}\label{def:ramification-morphism}
  A \emph{ramification morphism} is a map 
  $\rho :({\mathbb C}^{2},0)\rightarrow ({\mathbb C}^{2},0)$ such that there are local analytic
  coordinates $(u,v)$ and $(x,y)$ at $(0,0)\in ({\mathbb C}^{2},0)$ in which
  \begin{equation*}
    \left(x(u,v),y(u,v)\right)=\rho (u,v)=(u^{r},v)
  \end{equation*}
  for some positive integer $r$, called the \emph{ramification order}. These
  coordinates will be called \emph{adapted to $\rho $}.
\end{definition}
Notice that we allow different changes of  coordinates in each side of the
map. If $\mathcal{F}$ is a germ of foliation 
in $({\mathbb C}^{2},0)$, and $\rho^{\star}\mathcal{F}$ is its pull-back by 
$\rho $, then any separatrix of $\rho^{\star}\mathcal{F}$ is mapped by 
$\rho $ into a separatrix of $\mathcal{F}$
(although, obviously, this assignment needs  not be bijective).\par
All the arguments will be simplified if we introduce some notation: let
$\pi:\mathcal{X}\rightarrow ({\mathbb C}^{2},0)$  be a sequence of 
blowing-ups $\pi
=\pi _{k}\circ\dots\circ \pi_{1}$ whose respective exceptional
lines are $E_{i}$. We say that $\pi $ is a
\emph{regular tree of blowing-ups} if  all
the centers are regular points of the exceptional
divisor. The sequence $\pi $ is a
\emph{chain of blowing-ups} if  the center $P_{i+1}$ of $\pi _{i+1}$ belongs to
$E_{i}$. Finally, a regular tree which is a chain will be called a \emph{string of
  blowing-ups}. If $\pi $ is a regular tree, we say that an irreducible
component $E_{i}$ of the exceptional divisor $E$ is \emph{a son of} component
$E_{j}$ (and that $E_{j}$ is the father of $E_{i}$) if $E_{i}\cap E_{j}\neq
\emptyset$ and $i>j$.\par
Given $\pi $ as above and a ramification morphism $\rho $, we say that $\rho $
is \emph{transversal to $\pi $} if there are
coordinates $(u,v)$ and $(x,y)$ adapted to $\rho $ such that none of the centers
of $\pi_{i} $ is the infinitely near point given by $T_{(0,0)}(x=0)$. \par
Recall \cite{Hironaka} that given a rational function $g/h=f:\mathcal{Y}\rightarrow
{\mathbb C}$ where $\mathcal{Y}$ is an analytic surface, a point $P\in
\mathcal{Y}$ is an \emph{indetermination point of $f$} if  $g(P)=h(P)=0$. It is
well-known that there is a finite sequence of point blowing-ups
$\eta: \tilde{\mathcal{Y}}\rightarrow \mathcal{Y}$ such that $f\circ\eta$ is
well-defined everywhere; this $\eta$ is the \emph{elimination of
  indeterminations of $f$}. Obviously, this construction can be applied to maps
from $\mathcal{Y}$ to any analytic space (see \cite{Hironaka} for the details).\par
As our aim is to prove the Separatrix Theorem, we are going to assume from now
on that all our foliations are \emph{non-dicritical}.\par
Our main result is the following
\begin{theorem}
  \label{the:1}
  Let $\pi
  :\mathcal{X}\rightarrow ({\mathbb C}^{2},0)$ be the minimal reduction of
  singularities of $\mathcal{F}$. 
  There exists a ramification
  $\rho $ such that if $\tilde{\pi }:\tilde{\mathcal{X}}\rightarrow ({\mathbb
  C}^{2},0)$ is the elimination of indeterminations of $\pi ^{-1}\circ \rho
  :({\mathbb C}^{2},0)\rightarrow \mathcal{X}$,  then $\tilde{\pi }$ is a
  regular tree and all the singularities of $\tilde{\pi }^{\star}\rho
  ^{\star}\omega $ are simple.
\end{theorem}
As an obvious consequence, we get:
\begin{corollary}
  \label{cor:1}
  There is a ramification morphism $\rho $ such that the reduction of
  singularities of $\rho ^{\star}\omega $ is a regular tree.
\end{corollary}
Before giving the proof of Theorem \ref{the:1}, it is convenient to introduce
some more notation:
\begin{definition}\label{def:1}
  Let $E=\cup_{i=1}^{n} E_{i}$ be the exceptional divisor of $\pi $. We say that the
  irreducible component $E_{i}$ corresponding to $\pi _{i}$ is $\mathcal{F}$\emph{-terminal}
  if there is no $j>i$ such that the center $P_{j}$ of $\pi _{j}$ belongs to
  $E_{i}$.
\end{definition}
A germ of analytic curve $\gamma :({\mathbb C},0)\rightarrow ({\mathbb
  C}^{2},0)$ has the \emph{type of singularity} of $E_{i}$ if its strict
transform $\pi ^{\star}\gamma $ meets $E$ transversely at $E_{i}$ (as a consequence, $\pi
^{\star}\gamma $ is also regular). The germ $\gamma $ is $\mathcal{F}-$terminal
if it has the type of singularity of an $\mathcal{F}-$terminal divisor. 
\begin{remark}\label{rem:1}
  It is clear that 
  the reduction of singularities of a foliation $\mathcal{F}$ is a regular tree if
  and only if any  $\mathcal{F}-terminal$ curve is regular.  
\end{remark}
We shall make constant use of the following result, whose proof is an easy
computation:
\begin{lemma}\label{lem:1}
  Let $\pi $ be a sequence of blowing-ups and $\rho $ a ramification of order
  $r$ transversal to $\pi $. Call $\tilde{\pi }$ to the elimination of singularities of $\pi
  ^{-1}\circ \rho $. Then:
  \begin{itemize}
  \item[a)] If $\pi $ is a string of $n$ blowing-ups then $\tilde{\pi }$ is a
    string of $nr$ blowing-ups.
  \item[b)] As a consequence, if $\pi $ is a tree with $b$ branches of lengths
    $l_{1},\dots l_{b}$, then $\tilde{\pi }$ is a tree with $b$ branches of
    lengths $rl_{1},\dots, rl_{b}$.
  \end{itemize}
\end{lemma}
For the reader's sake, we give the local expression of $\tilde{\rho }=\pi
^{-1}\circ\rho\circ\tilde{\pi }$ when $\pi $ is a string. Fix systems of
coordinates $(u,v)$ and $(x,y)$ in $({\mathbb C}^{2},0)$ adapted to $\rho $. 
Cover $\mathcal{X}$ and $\tilde{\mathcal{X}}$
as follows:
\begin{equation}
  \label{eq:blowing-ups-charts}
  \left\{
    \begin{array}{l}
      \mathcal{X}_{n}=\left(\overline{U}_{1}\cup\dots\cup\overline{U}_{n}\right)\cup U_{n} \\
      \tilde{\mathcal{X}}=\left(\overline{U}_{11}\cup\dots\cup\overline{U}_{1r}\right)\cup\dots\cup\left(\overline{U}_{n1}
      \cup\dots\cup\overline{U}_{nr}\right)\cup
      U_{nr}\\           
    \end{array}
  \right.
\end{equation}
where each $\overline{U}_{i}$ is the standard chart covering \emph{the point at
infinity} of the exceptional divisor $E_{i}$ of $\pi_{i}$, and $U_{n}$ covers
the origin of $E_{n}$. The open set $\overline{U}_{ij}$ covers the point at
infinity of $E_{r\dot (i-1)+ j}$ and $U_{nr}$ covers the origin of
$E_{nr}$. These charts have respective systems of coordinates
$(\overline{x}_{i},\overline{y}_{i})$ and $(x_{n},y_{n})$ for the
$\overline{U}_{i}$'s and $U_{n}$, and $(\overline{u}_{ij},\overline{v}_{ij})$
and $(u_{nr},v_{nr})$ for the $\overline{U}_{ij}$ and $U_{nr}$, which can be
chosen (in a standard way) so that
\begin{equation}
  \label{eq:blowing-ups-charts-coord}
  \left\{
    \begin{array}{l}
  \left(\overline{u}_{nr},\overline{v}_{nr}\right)=\left(\dfrac{1}{v_{nr}},u_{nr}v_{nr}\right)
  \text{ and }\\[10pt]
  \left(\overline{u}_{i'j'},\overline{v}_{i'j'}\right)=\left(\dfrac{1}{v_{ij}},u_{ij}v_{ij}^{2}\right)
  \text{ when }
  r (i-1)+j = r (i'-1)+j'+1
  \end{array}
  \right.
\end{equation}
and the analog conditions in $\mathcal{X}$. We have (assuming $\pi $ is a string):
\begin{proposition}\label{pro:1}
  Let $\tilde{E}=\left(E_{11}\cup\dots\cup E_{1r}\right)\cup\dots\cup\left(E_{n1}\cup\dots\cup E_{nr}\right)$ be the
  exceptional divisor of $\tilde{\pi }$ and $E=E_{1}\cup \dots\cup E_{n}$ that of $\pi
  $. Call $P_{i}$ to the origin of the chart $\overline{U}_{i}$ (the point ``at
  infinity''). Then:
  \begin{enumerate}
  \item For $i=1,\dots,n$ and $j=1,\dots,r-1$, $\tilde{\rho }(E_{ij})=P_{i}$.
  \item For $i=1,\dots,n$, $\tilde{\rho }(E_{ir})=E_{i}$.
  \end{enumerate}
  Moreover, given $i$ and $j$ as above, The local expression of $\tilde{\rho }$ in
  $\overline{U}_{ij}$ and $\overline{U}_{i}$ is:
  \begin{equation}\label{eq:local-expression-for-strings}
    (\overline{x}_{i},\overline{y}_{i})=\tilde{\rho }(\overline{u}_{ij},\overline{v}_{ij})=
    \left(\overline{u}_{ij}^{r-j+1}\overline{v}_{ij}^{r-j}, \overline{u}_{ij}^{j-1}\overline{v}_{ij}^{j}\right)
  \end{equation}
\end{proposition}
Combining Proposition \ref{pro:1} with Lemma \ref{lem:1} and taking into account that the
transform of a simple singularity by a map of the form
(\ref{eq:local-expression-for-strings}) is also a simple singularity, one infers:
\begin{corollary}\label{cor:2}
  If $\rho $ is such that the reduction of singularities of $\rho
  ^{\star}\mathcal{F}$ is a regular tree, then for any other ramification
  $\sigma $,
  the reduction of singularities of $(\rho \circ \sigma )^{\star} \mathcal{F}$ is also a
  regular tree.
\end{corollary}
%%% \begin{remark}\label{rem:2}
%%%   If the reduction of singularities $\pi $ of $\mathcal{F}$ is regular tree, and
%%%   $\tilde{\mathcal{F}}_{P}$ is the germ at $P\in E$ of the strict transform of
%%%   $\mathcal{F}$ \emph{at any stage of $\pi $}, then the reduction of
%%%   singularities of $\tilde{\mathcal{F}}_{P}$ is, obviously, a regular tree.
%%% \end{remark}
\def\proofname{Proof of Theorem \ref{the:1}}
\begin{proof}
Let $T=F_{1}\cup \dots\cup F_{p}$ be the
union of all the terminal divisors of $\mathcal{F}$, and let $\gamma
_{1},\dots,\gamma _{p}$ be a corresponding family of $\mathcal{F}-$terminal
curves. We can assume, by Lemma \ref{lem:1} that all the $\gamma _{i}$ are
singular, and that $\pi $ is given by the sequence
$\pi =\pi _{p}\circ \dots \circ \pi _{1}$, where $\pi _{i}$ is (the remaining
part of) the sequence of blowing-ups leading to $T_{i}$. Again by Lemma
\ref{lem:1} and Remark \ref{rem:1}, it is enough to show that for any $i$, there
is a ramification morphism $\rho _{i}$ such that $\rho_{i}^{\star}\gamma _{i}$ is a
union of regular curves, for $\rho =\rho _{p}\circ\dots\circ\rho _{1}$ would
satisfy the thesis. Hence, we need only prove that given an irreducible singular curve
$\gamma$, there is a ramification $\rho $ such that 
\begin{enumerate}
\item The curve $\rho ^{\star}\gamma$ is a
  union of non-singular branches.
\item The elimination of indeterminations of $\pi_{1}
  ^{-1}\circ \rho$, say 
  $\tilde{\pi }: \tilde{\mathcal{X}}\rightarrow ({\mathbb C}^{2},0)$, is a
  regular tree.
\item If $Q\in \tilde{\mathcal{X}}$ is
  such that $\tilde{\rho }(Q)$ is simple for $\pi ^{\star}\omega $, then $Q$ is
  simple for $\tilde{\rho}^{\star}\pi _{i}^{\star}\omega $.
\end{enumerate}
If $\gamma$ has multiplicity $m$, then property $(1)$ holds obviously for any
ramification of order $r$ multiple of $m$.\par
Let $\Gamma $ be a curve with only one Puiseux exponent, having maximal
contact with $\gamma $, so that they share all the infinitely near points of the
reduction of singularities $\pi _{\Gamma }$ of $\Gamma $. Let $T_{\Gamma }$ be
the last exceptional divisor of $\pi _{\Gamma }$. Reasoning by induction on the
number of Puiseux exponents of $\gamma $ and using Lemma \ref{lem:1} again, we
%%% If we prove that there is
%%% $\rho_{\Gamma } $ holding $(1)$, $(2)$ and $(3)$ and such that the corresponding
%%% $\tilde{\rho}_{\Gamma }$ is a ramification morphism at any point, then by
%%% induction on the number of Puiseux exponents we finish. That is,
may assume
that $\gamma $ has a single Puiseux exponent, $m/n$ with $(m,n)=1$.\par
Let then $\rho $ be a ramification morphism of order $r$ transversal to $\gamma $ and write
the continuous fraction expansion of $m/n=[a_{0},a_{1},\dots,a_{s}]$. Let $\pi _{0}$ be the
string following the maximal contact sequence of $\gamma $, $\pi ^{0}$ the
remaining blowing-ups and consider the diagram
\begin{equation} 
  \label{eq:cd1} 
  \begin{split} 
    \xymatrix{ 
      & \mathcal{X}\ar[d]^{\pi^{0}}\\ 
      \tilde{\mathcal{X}}\ar@{-->}[ur]^{\tilde{\rho }}\ar[r]^{\tilde{\rho}_{0} 
      }\ar[d]^{\tilde{\pi}_{0}} &  
      \mathcal{X}_{0}\ar[d]^{\pi_{0}}\\ 
      \mathcal{B}_{2}\ar[r]^\rho & \mathcal{B}_{2} 
    } 
  \end{split} 
\end{equation}
wihere $\tilde{\mathcal{X}}$ is the elimination of singularities of $\pi
_{0}^{-1}\circ \rho $. We use the
notation of Proposition \ref{pro:1} for $\pi _{0}$.
We want to prove that for some $r$, the map $\tilde{\rho }_{0}$ lifts
holomorphically to $\mathcal{X}$ (in other words, the \emph{rational} map
$\tilde{\rho }$ is actually holomorphic).
The sequence $\pi ^{0}$ starts by blowing-up $E_{n}\cap E_{n-1}$, so that we need
only check that $\tilde{\rho }_{0}$ \emph{restricted to\/}
$E_{n-1r}\cup E_{n1}\cup\dots\cup E_{nr-1}$ lifts to $\mathcal{X}$. Consider the matrices
\begin{equation*}
  A=
  \begin{pmatrix}
    1 & -1\\
    0 & 1
  \end{pmatrix}\,\,\text{ and }\,\,
  B=
  \begin{pmatrix}
    1 & 0\\
    -1 & 1
  \end{pmatrix},    
\end{equation*}
and, for any $i=a_{1}+\dots + a_{l}+s$ let $\varphi _{i}(A,B)=B^{s}A^{a_{l}}\cdots
B^{a_{2}}A^{a_{1}}$ (assuming $l$ odd, and the corresponding product for $l$ even).
The equations of $\tilde{\rho }$ at $U_{nj}$ are
\begin{equation}
  \label{eq:4}
  \tilde{\rho }(u_{j},v_{j})=\left(u_{j}^{e_{1}}v_{j}^{e_{2}},u_{j}^{f_{1}}v_{j}^{f_{2}}\right)
\end{equation}
where
\begin{equation*}
  \begin{pmatrix}
    e_{1} & e_{2}\\
    f_{1} & f_{2}
  \end{pmatrix}=
  \varphi_{i} (A,B)\cdot
  \begin{pmatrix}
    r - j + 1 & r - j\\
    j - 1 & j
  \end{pmatrix}
\end{equation*}
for some $0\leq i \leq a_{1}+\dots + a_{s}$ which depends on the point of $\mathcal{X}$ we are looking
at. From this, one sees that $\tilde{\rho }$ is holomorphic in $U_{nj}$
if and only if $e_{1}e_{2}\geq 0$ and $f_{1}f_{2}\geq 0$ (i.e. if the exponents
in each component do not have different signs). We are seeking for an $r$ for
which all the possible pairs $(e_{1},e_{2})$ and
$(f_{1},f_{2})$ share this property. Another inductive argument shows that for
all $i$, there exist
integer numbers $c_{i},d_{i},m_{i}$ such that
\begin{equation*}
  \begin{pmatrix}
    e_{1} & e_{2}\\
    f_{1} & f_{2}
  \end{pmatrix}=
  \begin{pmatrix}
    r - c_{i}(m_{i}-1) & r - c_{i}m_{i}\\
    d_{i} (m_{i}-1) & d_{i}m_{i}
  \end{pmatrix}.
\end{equation*}
Hence, if $r$ is a multiple of $c_{i}d_{i}$ for all $0\leq i \leq a_{1}+\dots
+a_{s}$, then all the above maps
are well-defined and $\tilde{\rho }$ is holomorphic. As the number of $i$'s is
finite, this $r$ exists.
\end{proof}
\def\proofname{Proof}

%
%%% Local Variables: ***
%%% mode:tex ***
%%% TeX-master:"00father.ltx" ***
%%% comment-column:0 ***
%%% comment-start: "%%% "  ***
%%% comment-end:""
%%% End:

\section{Recalling Residues}
Given a germ of holomorphic foliation $\mathcal{F}$ in $({\mathbb C}^{2},0) $,
and a regular separatrix $E$ passing through $Q=(0,0)$, the \emph{Camacho-Sad} index of ${\mathcal F} $
at $Q$ along $E$ is defined as
$$
\text{ I}_{Q}(\mathcal{F},S)=-\text{ Res}_{0}\left(\frac{a(x,0)}{b(x,0)}dx\right).
$$
where we are assuming that ${\mathcal F} $ is defined near $Q$ by $\omega
=ya(x,y)dx + b(x,y)dy$ and $E\equiv (y=0)$.\par
One key result for any (known) proof of the Separatrix Theorem is the following
classical result (see \cite{Camacho-Sad-Coloquio} for a complete proof and
references to the original sources):
\begin{proposition}\label{pro:simple-separatrix}
  If $Q$ is a simple singularity of ${\mathcal F} $, $S$ is a smooth
  separatrix of ${\mathcal F} $ through $Q$ with $\text{ I}_{Q}(\mathcal{F},
  S)\neq 0$ then there is another (smooth) separatrix $T$ passing through $Q$.
\end{proposition}
Let now $\pi:{\mathcal X} \rightarrow ({\mathbb C}^{2},0)$ be the blowing-up of $(0,0)$ and
call $E$ to the exceptional divisor. Then
\begin{theorem}[\cite{Camacho-Sad}]\label{the:residue--1}
  If $E$ is invariant for the pull-back $\pi^{\star}\mathcal{F}$ of
  $\mathcal{F}$ and $P_{1},\dots,P_{r}$ are the singular points of
  $\pi^{\star}\mathcal{F}$ in $E$, then
\begin{equation}
    \label{eq:residue-blowup}
    \sum_{i=1}^{r}I_{P_{i}}(\mathcal{F},S)=-1
\end{equation}
Moreover, if $S$ is invariant for $\mathcal{F}$, 
$\pi^{\star}S$ is its pull-back by $\pi$ and $Q=E\cap S$, then
\begin{equation}
  \label{eq:residue-blowup2}
  I_{Q}(\pi^{\star}\mathcal{F},\pi^{\star}S)=I_{(0,0)}(\mathcal{F},S)-1.
\end{equation}
\end{theorem}
Finally, it is easily verified that if $S$ and $T$ are both smooth separatrices
of ${\mathcal F} $ and $Q=S\cap T$ is a simple singularity of ${\mathcal F} $,
then
\begin{equation}
  \label{eq:inverse-residue}
  \text{ I}_{Q}(\mathcal{F},S)=\frac{1}{\text{ I}_{Q}({\mathcal F} ,T)}
\end{equation}
whenever any of both indices is \emph{non-zero}, or using the standard
terminology, when $Q$ is not a saddle-node singularity of $\mathcal{F}$.

%
%%% Local Variables: 
%%% mode:tex 
%%% TeX-master:"00father.ltx" 
%%% comment-column:0 
%%% comment-start: "%%% "  
%%% comment-end:""
%%% End:

\section{Proof of the Separatrix Theorem}
We have now all the machinery needed to give a straightforward proof of the
Separatrix Theorem:
\begin{theorem}[\cite{Camacho-Sad}]
  Given a germ of holomorphic foliation ${\mathcal F} $ at $({\mathbb C}^{2},0)
  $, there is a separatrix for ${\mathcal F} $ passing through $(0,0)$.  
\end{theorem}
\begin{proof}
As stated in Section \ref{sec:ramifications}, ramifications send separatrices of
$\rho ^{\star}{\mathcal F} $ into separatrices of ${\mathcal F} $ so that, by
Corollary \ref{cor:1}, we need only prove the result for foliations whose
reduction of singularities is a regular tree.\par
Assume then, that the reduction of singularities $\pi :{\mathcal X} \rightarrow
({\mathbb C}^{2},0) $ of ${\mathcal F} $ is a tree and let $E=E_{1}\cup\dots
\cup E_{n}$ be the exceptional divisor. By Proposition
\ref{pro:simple-separatrix}, we only need to show that there is a point $Q\in
E_{i}$ for some $i$ with $Q\not\in E_{j}$ for $i\neq j$ such that $\text{
  I}_{Q}({\mathcal F} , E_{i})\neq 0$.\par
Assume, by contradiction, that there is no such $Q$. That is, for any $Q$ in the
regular part of $E$, the corresponding index is $0$. We have:
\textbf{Assertion:\/} Given an irreducible component $F$ of $E$, if
$F_{1},\dots, F_{s}$ are its sons ($s$ may be $0$), then
\begin{equation}
  \label{eq:sons}
  \sum_{Q\in F}\text{ I}_{Q}(\mathcal{F}, F)=-s-1,
\end{equation}
and all the terms in the sum are \emph{rational numbers}.
This is (easily) proved by induction on the maximal length of a branch starting from $F$,
using Theorem \ref{the:residue--1} and our assumption that all the
residues at the regular part of $F$ are $0$. A direct consequence of this is
that if $F'$ is the father of $F$ and $P=F\cap F'$, then
\begin{equation}
  \label{eq:inequality}
  \text{ I}_{Q}({\mathcal F} ,F')\geq -1
\end{equation}
(and the above index is a rational number, in fact).\par
But now, the divisor $E_{1}$ appearing in the first blowing-up has no father, so
that all the singularities belonging to it either are the crossing with a son or
have zero index along $E_{1}$. From the Assertion and Equation
(\ref{eq:inequality}), if $F_{1},\dots,F_{s}$ are the sons of $E_{1}$ and
$Q_{i}=E_{1}\cap F_{i}$, we must have
\begin{equation*}
  \begin{split}
      -r-1 &\stackrel{(\ref{eq:sons})}{=} \sum \text{ I}_{Q}({\mathcal F},
      E_{1})\\
      &{=}\sum_{i=1}^{s}\text{ I}_{Q_{i}}({\mathcal F}
      ,E_{1})\\
      &\stackrel{(\ref{eq:inequality})}{\geq} \sum_{i=1}^{s}-1 = -r
  \end{split}
\end{equation*}
which gives the desired contradiction.
\end{proof}

%
%%% Local Variables: 
%%% mode:tex 
%%% TeX-master:"00father.ltx" 
%%% comment-column:0 
%%% comment-start: "%%% "  
%%% comment-end:"" 
%%% End: 

\providecommand{\bysame}{\leavevmode\hbox to3em{\hrulefill}\thinspace}

%\bibliographystyle{amsplain}
%\bibliography{../biblioams}
\end{document}